\documentclass{article}

\usepackage{amsmath,amsfonts,amsthm,units}
\usepackage{hyperref}                     
\hypersetup{colorlinks=true}              

\theoremstyle{plain} \newtheorem{prop}{Proposition}[section]
\theoremstyle{plain} \newtheorem{lemma}[prop]{Lemma}        
\theoremstyle{plain}     
\theoremstyle{plain} \newtheorem{defin}[prop]{Definition}   
\theoremstyle{plain} \newtheorem{thm}[prop]{Theorem}        
\theoremstyle{plain} \newtheorem{remark}[prop]{Remark}      
\theoremstyle{plain}     
\theoremstyle{plain}  
\theoremstyle{definition}    
\theoremstyle{definition}  

\newcommand{\sgn}{\operatorname{sgn}}                             
\newcommand{\Reals}{{\mathbb R}}
\newcommand{\reals}{{\mathbb R}}

\newcommand{\One}[2]{\mathbf{1}_{\left[ #1,#2 \right]}}           
\newcommand{\tends}[2]{\stackrel{#1\rightarrow #2}{\longrightarrow}}
\newcommand{\norm}[1]{\left\Vert#1\right\Vert}                    
\newcommand{\abs}[1]{\left\vert#1\right\vert}                     
\newcommand{\ud}{\,\mathrm{d}}

\newcommand{\goes}{\stackrel{n\rightarrow\infty}{\longrightarrow}}

\newcommand{\rkhs}{\mathcal{H}}
\newcommand{\bfrac}[2]{\frac{\mbox{\small{$#1$}}}{\mbox{\small{$#2$}}}}
\newcommand{\E}{\mathbb{E}}

\newcommand\blfootnote[1]{%
  \begingroup
  \renewcommand\thefootnote{}\footnote{#1}%
  \addtocounter{footnote}{-1}%
  \endgroup
}

\begin{document}


\title{Covariance of stochastic integrals with respect to fractional Brownian motion}%
\author{Yoha\"{i} Maayan
\footnote{Department of Mathematics, Technion, Israel Institute of Technology, 32000 Haifa, Israel}
 \and  
Eddy Mayer-Wolf\footnotemark[1]}


\maketitle

\begin{abstract}
We find an explicit expression for the
cross-covariance between stochastic integral processes
with respect to a $d$-dimensional fractional Brownian motion (fBm) $B_t$ with Hurst
parameter $H\!>\!\nicefrac{1}{2}$, where the integrands are vector fields applied to $B_t$. 
It provides, for example, a direct alternative proof of Y. Hu and D. Nualart's result that the stochastic integral component in the fractional Bessel process decomposition is not itself a fractional Brownian motion.
\end{abstract}





\section{Introduction}
Fractional Brownian motion is a family of zero mean stationary Gaussian processes $B_t\!=\!B_t^H$ indexed by
$H\!\!\in\!\!\left(0,1\right)$ which was mathematically introduced by B.B Mandelbrot
and J.W. Van Ness in~\cite{Mandelbrot}\ (cf.~\cite{Kolmogorov} as well).
It generalizes Brownian motion ($H\!=\!\nicefrac{1}{2}$) in that $EB_t^2=t^{2H}$, and can be used to model various phenomena, in finance as well as in other fields. This is primarily due to the
fact that its self-similarity depends on the parameter $H$, 
which allows for phenomena exhibiting different kinds of self-similarity to be
modeled by fractional Brownian motion with an appropriate $H$.

Since fractional Brownian motion is not a semimartingale (unless
$H\!=\!\nicefrac{1}{2}$), the ordinary stochastic
calculus for semimartingales (such as the It\^{o} integral) does
not apply. Instead, there are several approaches for defining a
stochastic integral with respect to fractional Brownian motion.
The divergence integral is one possible approach, the one discussed in this paper, using the
Malliavin divergence operator as the basis for integration, a
survey of which can be found in D. Nualart's book~\cite{Nualart}. One other approach for example was developed by Z\"{a}hle in~\cite{Zahle}, which involves a
pathwise definition of the stochastic integral. This requires a
generalization of the Young-Stieltjes integral, introduced in
the same paper.

Given a suitable process $u$, the divergence integral
$t\mapsto\int_0^tu_\tau\ud B_\tau$ yields a new process $X_t$. While the
general theory of Malliavin calculus provides an abstract formula for
the covariance function of $X_t$ (c.f.~\cite{Nualart}), in many concrete cases
it is not straightforward to write it explicitly.
In this paper such an expression is provided for the correlation between $\int_0^tF\left(B_\tau\right)\ud B_\tau$
and $\int_0^sG\left(B_\sigma\right)\ud B_\sigma$. 

In~\cite{Hu_Nualart} Y. Hu and D. Nualart showed that if $H\!\ne\!\nicefrac{1}{2}$ the process
$X_t\!=\!\int_0^t \operatorname{sgn}(B_s)\ud Bs$ is not a fractional Brownian motion (and similarly in the multidimensional case). This is different than the case $H=\frac12$, where these processes (in any dimension) are also Brownian motions.
In the absence of a formula for the covariance of this process, a detailed analysis of its chaos
expansion was necessary to reach that conclusion, and this was one of the motivations to obtain
such a formula which, moreover, would indeed have to accommodate non-smooth functions $F$ such as $\operatorname{sgn}(x)$. This is addressed in Section~\ref{Bessel}.

Section~\ref{sec:fBm} includes some preliminaries and an
auxiliary result (Proposition~\ref{lemma:kernels}).
The main covariance formula is presented in
Section~\ref{gen_res}. It is first stated for relatively regular $F$ and $G$ as Theorem~\ref{theorem:special_case_of_main}, and then in its full generality as Theorem~\ref{theorem:main}.


\section{Preliminaries: Fractional Brownian Motion}  \label{sec:fBm}
The following introduction to fractional Brownian motion and its analysis is taken mostly from Chapter 5 in~\cite{Nualart}.
\begin{defin}
A fractional Brownian motion with Hurst parameter $H\!\in\!(0,1)$ is a centered Gaussian process $\left\{B_t,\,t\in\left[0,T\right]\right\}$ defined on a complete probability space $(\Omega,{\mathcal F},\mathbb{P})$ with covariance function
\begin{equation}\label{eq:fBmcov}
R_H\left(t,s\right)=\frac{1}{2}\left(t^{2H}+s^{2H}-\abs{t\!-\!s}^{2H}\right).
\end{equation}
A vector process $B_t\!=\!\left(B^1_t,\ldots,B^d_t,\,t\in [0,T]\right)$ with independent fractional Brownian motion components will be referred to as a $d$-dimensional fractional Brownian motion.
\end{defin}
\noindent
The parameter $H$ will not be explicitly indicated in the notation $B_t$.
It follows that
\begin{equation}
E\abs{B_t-B_s}^2=\abs{t\!-\!s}^{2H}
\end{equation}
and, by Kolmogorov's continuity criterion, we may assume that
$\left\{B_t\mid t\in\left[0,T\right]\right\}$ has
$\alpha$-H\"{o}lder continuous trajectories for any $\alpha <H$.
When $H=\nicefrac{1}{2}$, $B_t$ is a standard Brownian
motion.

From this point on, it will be assumed that $H>\frac{1}{2}$. The
discussion takes a different turn when $H<\frac{1}{2}$ (more about that in Section~\ref{sec:concluding}).
It will be convenient to write
\begin{equation}\label{eq:integral.cov}
 R_H\left(t,s\right)=\alpha_{H}\int_0^t\int_0^s\abs{\tau\!-\!\sigma}^{2H-2}\ud \tau\ud \sigma,
\end{equation}
 with
\begin{equation} \label{alpha}
 \alpha_H=H\left(2H-1\right),
\end{equation}
as can be easily verified.
 Denote by $\mathcal{E}$ the space of step functions (that is, the space spanned by all
 indicators of subintervals of $\left[0,T\right]$), and on it define
\begin{equation}\label{eq:integral.inner.product}
\left<f,g\right>=\alpha_H\iint_{\left[0,T\right]^2}f\left(\tau\right)g\left(\sigma\right)\abs{\tau-\sigma}^{2H-2}\ud
\tau\ud \sigma\ ,
 \hspace{.7cm} f,g\in\mathcal{E}.
\end{equation}
In particular,
$\left<\mathbf{1}_{\left[0,t\right]},\mathbf{1}_{\left[0,s\right]}\right>=R_H\left(t,s\right)$.
The completion $\mathcal{H}$ of $\mathcal{E}$ with respect to
$\left<\cdot,\cdot\right>$ is the ``deterministic integrands space'' associated with $B_t$ (or with $R_H(s,t)$).

\noindent
The linear isometry
$\mathbf{1}_{\left[0,t\right]}\!\in\!\mathcal{E}\longmapsto \delta\mathbf{1}_{\left[0,t\right]}:=B_t\!\in\!L^2(\mathbb{P})$
 extends to $\mathcal{H}$.
The Malliavin derivative $D$ is then defined as an unbounded operator from $L^2\left(\mathbb{P}\right)$ to
$L^2\left(\mathbb{P};\mathcal{H}\right)$ by setting $D\delta h=h$ and
extending it, by the chain rule
\begin{equation*}
 f\left(\delta h_1,\ldots,\delta h_n\right)\longmapsto
   \sum_{i=1}^n\partial_if\left(\delta h_1,\ldots,\delta h_n\right)h_i
\end{equation*}
and closure, to its domain $\mathbb{D}^{1,2}$ (the same can be
done for any $p\geq 1$ instead of $p=2$).\ The Malliavin derivative
of $F\in L^2\left(\mathbb{P};K\right)$, for a separable Hilbert space $K$,
is defined similarly as an element of
$L^2\!\left(\mathbb{P};\mathcal{H}\otimes K\right)$, its domain denoted
$\mathbb{D}^{1,2}\!\left(K\right)$ accordingly. We shall identify $\mathcal{H}\otimes K$ with the space of Hilbert-Schmidt operators from $K$ to $\mathcal{H}$.

The dual operator\ $\delta u$\ of\ \ $D\!:\!L^2(\mathbb{P})\!\to\!L^2(\mathbb{P};\mathcal{H})$
\,, commonly referred to as the stochastic integral $\int_0^Tu_t\ud B_t$, satisfies by definition
\begin{equation*}
E\left<DF,u\right>=EF\delta u\qquad\forall F\in\mathbb{D}^{1,2}.
\end{equation*}
Its domain $\operatorname{Dom}\left(\delta\right)$ contains
$\mathbb{D}^{1,2}\left(\mathcal{H}\right)$ and moreover
\begin{equation}\label{absfor}
  E\delta u\delta v = E\left<u,v\right> + E\operatorname{trace}\left(Du Dv\right) ,
    \hspace{0.7cm}u,v\in\mathbb{D}^{1,2}\left(\mathcal{H}\right)
\end{equation}
(cf. \cite[Proposition 1.3.1]{Nualart}).\ \
 If $u\!\equiv\!h\in\mathcal{H}$
the notation $\delta u$ is consistent with the previous usage of $\delta h$,
and the term ``integral'' reflects the fact that when
$H\!=\!\nicefrac{1}{2}$, $\delta u$ coincides with the Skorohod
integral of $u_t$ with respect to Brownian motion
.
 \begin{remark}\label{multi_dim_malliavin}
	For the sake of simplicity, ${\mathcal H}$, $D$ and $\delta$ were defined only for a scalar fractional Brownian
    motion. The integrands space\ \ ${\mathcal H}_d$ associated with a $d$-dimensional fractional Brownian motion
    is the Hilbert direct sum of $d$ copies of the one-dimensional integrands spaces\ \ ${\mathcal H}$. For appropriate $F$ and $u$ (respectively scalar and ${\mathcal H}_d$-valued),\ \  $DF\!=\!\left(D^1F,\cdots, D^dF\right)$\,,\
    $\delta u\!=\!\delta^1u^1+\cdots+\delta^du^d$\ and\
	 $Du\,h\!=\!\left((D^1u^1)h^1,\ldots,(D^du^d)h^d\right)$\,,\ $h\!\in\!{\mathcal H}_d$.
 \end{remark}
Note that $\mathcal{H}$ contains not only proper functions but distributions as well. The function space
 \begin{equation}  \label{absH}
  \abs{\mathcal{H}}=\left\{f:\left[0,T\right]\rightarrow\mathbb{R}\ {\Big|}\ \iint_{\left[0,T\right]^2}
  \abs{f\left(t\right)}\abs{f\left(s\right)}\abs{t-s}^{2H-2}\ud s\ud t<\infty\right\}
 \end{equation}
is a (strict,~\cite{Pipiras&Taqqu}) dense subspace of $\mathcal{H}$,
on which the inner product is still given by \eqref{eq:integral.inner.product}.


\noindent
Similarly it will be convenient to single out the following subclass of elements of $\mathcal{H}\otimes\mathcal{H}$ whose inner product can be expressed explicitly.
\begin{defin}
 Let $K\in\mathcal{H}\!\otimes\!\mathcal{H}$. If\ $K(\mathcal{H})\subset\abs{\mathcal{H}}$\ and there exists $k(s,t)\!\in\!L^\infty([0,T]^2)$ such that for almost all
 $t\!\in\![0,T]$
  \begin{equation} \label{kernel_representation}
   Kh(t)=\left<k(\cdot\,,\,t)\,,\,h\right>_\mathcal{H},\hspace{1.5cm}h\!\in\!\mathcal{H},
  \end{equation}
  (in particular\ \ \ $\displaystyle{Kh(t)=
          \alpha_H\iint_{_{\mbox{\scriptsize $[0,T]^2$}}}\!\!\!\!\!k(s,t)\,h(u)\,|s\!-\!u|^{2H-2}\,\ud s \ud u}$\ \
              \ \ \ for all $h\!\in\!|\mathcal{H}|$)\\
   we shall say that $K$ is {\em represented by the kernel $k$}\ \ or simply {\em representable}.
\end{defin}
%
\begin{prop}\label{lemma:kernels}
Assume that $K_1,K_2\in\mathcal{H}\!\otimes\!\mathcal{H}$ are respectively represented by kernels $k_1, k_2\!\in\!L^\infty([0,T]^2)$
 . Then
   \begin{equation}\ \ \label{kernel_inner_product}
      \left<K_1,K_2\right>_{\mathcal{H}\otimes\mathcal{H}}
              =\alpha_H^2\int_{_{\mbox{\scriptsize $\left[0,T\right]^4$}}}\!\!\!\!\!\!\!\!\!k_1\left(s,t\right)
                                      k_2\left(u,v\right)\abs{s-u}^{2H-2}
                                              \abs{t-v}^{2H-2}\ud s\ud t\ud u\ud v.
   \end{equation}
Moreover, the trace-class operator $K\!=\!K_2K_1$ is represented by the kernel
 \begin{equation}\label{eq:kernel}
  k\left(s,t\right)
     =\alpha_H\iint_{_{\mbox{\scriptsize $[0,T]^2$}}}\!\!\!\!\!k_1\left(s,u\right),k_2\left(v,t\right)\,|u-v|^{2H-2}\,\ud u \ud v
 \end{equation}
  and its trace is given by
 \begin{equation}\label{eq:trace}
  \operatorname{trace}\left(K\right)
       =\alpha_H\iint_{_{\mbox{\scriptsize $[0,T]^2$}}}\!\!\!\!\!k\left(s,t\right)\abs{s-t}^{2H-2}\ud s\ud t.
 \end{equation}
\end{prop}
\begin{proof}
    When $K_i=h^i_1\otimes
  h^i_2$ where $h^i_j\in\abs{\mathcal{H}}$ and $i,j=1,2$,~\eqref{kernel_inner_product} is
  a simple calculation. This clearly carries over to sums of such operators. For general $k_1$ and $k_2$,
  let $\left\{k_i^n\right\}_{n=1}^{\infty}$ be a sequence of kernels of the type above which converges almost
  everywhere to $k_i$ and is uniformly bounded (this sequence exists since $k_i$ itself is bounded). Then each
  kernel defines a Hilbert Schmidt operator $K_i^n$, and~\eqref{kernel_inner_product} is known for $K_i^n$. It remains to
  prove that the sequence $\left\{K_i^n\right\}_{n=1}^{\infty}$ converges in the Hilbert
  Schmidt norm to $K_i$, and that the right hand side of~\eqref{kernel_inner_product} converges accordingly.

  To see that $K_i^n\tends{n}{\infty} K_i$, note that $\left\{K_i^n\right\}_{n=1}^{\infty}$ is a Cauchy sequence, as implied by
  dominated convergence with the already established formula in~\eqref{kernel_inner_product} for the Hilbert Schmidt norm $\norm{K_i^m-K_i^n}$.  For any
  $h\in\abs{\mathcal{H}}$,
  \begin{equation*}
   K_i^nh\left(t\right)
       =\alpha_H\iint_{\left[0,T\right]^2}h\left(s\right)k_i^n\left(\tau,t\right)\abs{s-\tau}^{2H-2}\ud s\ud
       \tau
  \end{equation*}
   which converges in $\mathcal{H}$ to $K_ih$ by dominated
   convergence. Thus $K_i^n\tends{n}{\infty} K_i$ and
   \begin{equation*}
    \lim_{n\rightarrow\infty}\left<K_1^n,K_2^n\right>=\left<K_1,K_2\right).
   \end{equation*}
  To see that
  \begin{align*}
   \lim_{n\to\infty}&\alpha_H\int_{_{\left[0,T\right]^4}}\!\!\!\!\!\!\!\!k_1^n\left(s_1,t_1\right)k_2^n\left(s_2,t_2\right)
                                                                  \abs{s_1-s_2}^{2H-2}\abs{t_1-t_2}^{2H-2}\!\!\!\ud s_1\!\ud s_2\!\ud t_1\!\ud t_2\\
   &=\alpha_H\int_{_{\left[0,T\right]^4}}\!\!\!\!\!\!\!\!k_1\left(s_1,t_1\right)k_2\left(s_2,t_2\right)
                                                  \abs{s_1-s_2}^{2H-2}\abs{t_1-t_2}^{2H-2}\!\!\!\ud s_1\!\ud s_2\!\ud t_1\!\ud t_2,
  \end{align*}
  note that the integrand on the left hand side is bounded by some
  constant multiple of the integrable function
  $\abs{s_1-s_2}^{2H-2}\abs{t_1-t_2}^{2H-2}$. An additional application of dominated convergence concludes the proof of~\eqref{kernel_inner_product}.
 For any $h\!\in\!\mathcal{H}$\ \ $K_1h(u)=\left<k_1(\cdot,u)\,,\,h\right>_\mathcal{H}$\ \ and thus
 \begin{eqnarray*}
  K_2(K_1h)(t)
   &=&\alpha_H\iint_{_{\mbox{\scriptsize $[0,T]^2$}}}
         k_2(v,t)\left<k_1(\cdot,u)h\right>_\mathcal{H}|u-v|^{2H-2}\,\ud u \ud v\\
   &=&\left<\alpha_H\iint_{_{\mbox{\scriptsize $[0,T]^2$}}}k_2(v,t)k_1(\cdot,u)|u-v|^{2H-2}\,
                                                              \ud u \ud v\,,\,h\right>_\mathcal{H}
 \end{eqnarray*}
 which proves the second assertion.

 \noindent
   For the trace, note that $K_1$'s adjoint $K_1^{\ast}$ is represented by $k_1^{\ast}(s,t)\!:=\!k_1(t,s)$.  From~\eqref{kernel_inner_product}
  \begin{eqnarray*}
   \operatorname{trace}\left(K\right)
          &\!\!=\!\!&\left<K_1^{\ast},K_2\right>=\alpha_H^2\iint\!\!\!\iint_{_{\mbox{\scriptsize $[0,T]^4$}}}
                  \!\!\!\!\!\!\!k_1(t,s)\,k_2(u,v)\,|s\!-\!u|^{2H-2}\,|t\!-\!v|^{2H-2}\,\ud s\ud t\ud u\ud v\\
          &\!\!=\!\!&\alpha_H\iint_{_{\mbox{\scriptsize $[0,T]^2$}}}\!\!\!\!\!\!\!k(t,v)\,|t\!-\!v|^{2H-2}\,\ud t \ud v\ .
  \end{eqnarray*}
This completes the proof.
 \end{proof}

\section{The Covariance Formula}\label{gen_res}
We recall from Section~\ref{sec:fBm} that a process $u$ in
$\mathbb{D}^{1,2}$ belongs to the domain of $\delta$ and that for
any $u,v\!\in\!\mathbb{D}^{1,2}$, formula~\eqref{absfor} holds,
namely
\begin{equation}\label{absfor:reg}
  E\int_0^T\!\!\!u_\tau\,\ud B_\tau\int_0^T\!\!\!v_\tau\,\ud B_\tau
        =E\left<u,v\right>+E\operatorname{trace}\left(Du Dv\right),
\end{equation}
 where $Du$ is viewed as a Hilbert-Schmidt operator--valued random variable.
 Our aim is to find a concrete expression for the
right hand side of~\eqref{absfor:reg} when the integrands are
respectively of the form
  \begin{equation}  \label{specialFG}
     u\!=\!F\!\left(B_\cdot\right)\One{0}{t}(\cdot)\ \ \ \ \text{and}\ \ \ \
     v\!=\!G\!\left(B_\cdot\right)\One{0}{s}(\cdot).
  \end{equation}
In fact, we will do this for the multi-dimensional case, in which case formula~\eqref{absfor:reg} amounts to (see Remark~\ref{multi_dim_malliavin})
 \begin{equation}\label{multi_dim_abs_for}
 E\int_0^T\!\!\!u_\tau\,\cdot\ud B_\tau\int_0^T\!\!\!v_\tau\,\cdot\ud B_\tau
        =\sum_{i=1}^dE\left<u_i,v_i\right>+\sum_{i,j=1}^dE\operatorname{trace}\left(D^iu_j D^jv_i\right).
 \end{equation}

 \bigskip

 \noindent
 Our main result if Theorem~\ref{theorem:main} below, of which we first state a particular case as Theorem~\ref{theorem:special_case_of_main} which makes easier reading.

Recall the constant $\alpha_H\!=\!H\left(2H\!-\!1\right)$ . Let
 \begin{equation}\label{gamma}
  	\gamma(t,s)=H\alpha_H\left((s\vee t)^{2H-1}-|t\!-\!s|^{2H-1}\right)\,
                \left((s\wedge t)^{2H-1}+|t\!-\!s|^{2H-1}\right)\ \ \ \ s,t\!>\!0
 \end{equation}
 which is positive, symmetric and for all $t\!>\!0$ satisfies\ $\gamma(t,s)\!\approx\!\alpha_H^2t^{4H-3}s$\ \ as $s\!\to\!0$. 
  \begin{thm}\label{generalization}\label{theorem:special_case_of_main}
  Let $B_t,\ t\!\ge\!0$\ be a fractional Brownian
  motion with $H\!>\!\nicefrac{1}{2}$ and let
  $F,G:\reals^d\!\to\!\reals^d$ be locally integrable functions with first order distributional derivatives that are functions, and such that for some constants $M,\beta>0$,
  \begin{equation}\label{eq:exponential_growth_assumption}
  \abs{DF\left(x\right)},\ \abs{DG\left(x\right)}\leq Me^{\beta\abs{x}} ,\ \ i,j=1,\ldots,d.
  \end{equation}
 Then for $t,s\!\in\![0,T]$\ the processes $F\left(B_{r}\right)\One{0}{t}$ and $G\left( B_{r}\right)\One{0}{s}$ are $\abs{\mathcal{H}}^d$-valued and belong to $\operatorname{Dom}\delta$,
\ and moreover
  \begin{align}\nonumber
    \E\!\!\left(\int_0^t\!\!\!F\left(B_r\right)\!\!\ud B_r\!\int_0^s\!\!\!G\left(B_r\right)\!\!\ud B_r\right)
    =\alpha_H\!\!\!\int_0^t\!\!\!\int_0^s\E &\left(F\left(B_{\tau}\right)\cdot G\left(B_{\sigma}\right)\right)\abs{\tau\!-\!\sigma}^{2H-2}
                                             \!\ud\sigma\!\ud\tau\\
                                             &\hspace{-3cm}+\sum_{i,j=1}^d\int_0^t\!\!\!\int_0^s\gamma\left(\tau,\sigma\right)\E\left(\partial_iF_j\left(B_{\tau}\right)\partial_jG_i\left(B_{\sigma}\right)\right)\ud\sigma\ud\tau. \label{form}
  \end{align}
\end{thm}
\medskip
To state this result in full generality (where the derivatives of $F$ and $G$ could be measures), denote
\begin{equation}\label{L}
L_{t,s}\left(x,y\right)=\int_{0}^{t}\int_0^{s}f_{\tau,\sigma}\left(x,y\right)\gamma\left(\tau,\sigma\right)\ud \tau\ud \sigma,\ \ x,y\in\reals^d,\ t,s\in\left[0,T\right]
\end{equation}
where $f_{t,s}$ is the joint density of $\left(B_t,B_s\right)$.
Then:
 \begin{thm}\label{theorem:main}
Let $B_t,\ t\!\ge\!0$\ be a fractional Brownian motion with $H\!>\!\nicefrac{1}{2}$ and assume
\begin{enumerate}
\item
$F,G:\reals^d\!\to\!\reals^d$ are measurable and have exponential growth at most: 
there are some $C,\beta\!>\!0$ such that
\begin{equation} \label{exponential_growth}
\abs{F(x)},\ \abs{G\left(x\right)}\le Ce^{\beta |x|}\ ,\ \ \ \ \ \ x\!\in\!\reals^d.
\end{equation}

\item\label{condition:continuity_of_convolution}
The first order distributional derivatives of $F$ and $G$ are measures and for $0\leq t,s\leq 1$, letting $H=\begin{pmatrix}
F\\
G
\end{pmatrix}$, $L_{t,s}\ast\left(DH\otimes DH\right)$ is well defined in a neighbourhood of the origin and continuous at the origin.

\end{enumerate}
 Then for $t,s\!\in\![0,T]$\ the processes $F\left(B_{r}\right)\One{0}{t}$ and $G\left( B_{r}\right)\One{0}{s}$ are $\abs{\mathcal{H}}^d$-valued and belong to $\operatorname{Dom}\delta$,
\ and moreover
  \begin{align}\nonumber
    \E\!\!\left(\int_0^t\!\!\!F\left(B_r\right)\!\!\ud B_r\!\int_0^s\!\!\!G\left(B_r\right)\!\!\ud B_r\right)
    =\alpha_H\!\!\!\int_0^t\!\!\!\int_0^s\E &\left(F\left(B_{\tau}\right)\cdot G\left(B_{\sigma}\right)\right)\abs{\tau\!-\!\sigma}^{2H-2}
                                             \!\ud\sigma\!\ud\tau\\
                                             &\hspace{-2cm}+\sum_{i,j=1}^d\int_{\reals^{2d}}L_{t,s}\left(x,y\right)\partial_iF_j\left(\ud x\right)\partial_jG_i\left(\ud y\right).\tag{\ref{form}'}\label{form:general}
\end{align}
\end{thm}
\bigskip

Before proving Theorem~\ref{theorem:main}, we show in Proposition~\ref{proposition:continuity_of_convolution} below that Theorem~\ref{theorem:special_case_of_main} is in fact a particular case of Theorem~\ref{theorem:main}. It should be noted that this is not simply an exercise in generalisation. Section~\ref{Bessel} deals with the natural case $F\left(x\right)=G\left(x\right)=\operatorname{sgn}\left(x\right)$ (and its multidimensional counterpart) which is not covered by the assumptions of Theorem~\ref{theorem:special_case_of_main}. In this case\ $F'=2\delta_0$ and Condition~\ref{condition:continuity_of_convolution} in Theorem~\ref{theorem:main} is really that $L\left(x,y\right)$ be continuous at the origin. Since
  \begin{equation} \label{finite.integral}
   \int_0^T\int _0^T \frac{\gamma\left(\tau,\sigma\right)}
   {\tau^H\sigma^H\sqrt{1-\rho^2\left(\tau,\sigma\right)}}\ud \tau\ud \sigma<\infty,
  \end{equation}
($\rho\left(\tau,\sigma\right)$ is the correlation between $B_\tau$ and $B_\sigma$) which follows from
a standard asymptotic analysis near the singularity lines, the continuity of $L\left(x,y\right)$ can now be deduced in $d=1$ from dominated convergence.
\begin{prop}\label{proposition:continuity_of_convolution}
If $p,q:\reals^d\to\reals$ have exponential growth at most: 
\begin{equation*}
\abs{p\left(x\right)},\abs{q\left(x\right)}\leq Ce^{\alpha\abs{x}}, 
\end{equation*}
then the convolution $L\ast\left(p\otimes q\right)$ exists and is everywhere continuous.
\end{prop}
\begin{proof}
First, for any $u,v\in\reals^d$,
\begin{align}\nonumber
\iint_{\reals^{2d}}L\left(x,y\right)\abs{p\left(u-x\right)}\abs{q\left(v-y\right)}\ud x\ud y&\leq C^2\iint_{\reals^{2d}}L\left(x,y\right)e^{\alpha\abs{u-x}+\abs{v-y}}\ud x\ud y\\
&\hspace{-3.7cm}=C^2\iint_{\left[0,1\right]^2}\gamma\left(s,t\right)E\left(e^{\alpha\abs{u-B_s}+\alpha\abs{v-B_t}}\right)\ud s\ud t.\label{absolute_integrability_for_convolution_kernel}
\end{align}
Since (thinking of the $\ell^1$-norm on $\reals^d$ for simplicity)
\begin{align*}
E\left(e^{\alpha\abs{u-B_s}+\alpha\abs{v-B_t}}\right)&\leq e^{\alpha\left(\abs{u}+\abs{v}\right)}\left(E\left(e^{\alpha\abs{B_s^1}+\alpha\abs{B_t^1}}\right)\right)^d\\
&\leq e^{\alpha\left(\abs{u}+\abs{v}\right)}\left(E\left(e^{\alpha\max_{0\leq t\leq 1}\abs{B_t^1}}\right)\right)^d<\infty
\end{align*}
(see for example~\cite[Theorem~3.2]{Adler_monograph}), the right-hand-side of~\eqref{absolute_integrability_for_convolution_kernel} is finite. Thus $g:=L\ast\left(p\otimes q\right)$ is well defined. Note that
\begin{equation}
g\left(u,v\right)=\iint_{\left[0,1\right]^2}\gamma\left(s,t\right)E\left(p\left(u-B_s\right)q\left(v-B_t\right)\right)\ud s\ud t.
\end{equation}
Since the estimates above provide a domination for the $s$-$t$ integration, it will suffice to show that 
\[ E\left(p\left(u-B_s\right)q\left(v-B_t\right)\right)\underset{\left(u,v\right)\to\left(u_0,v_0\right)}{\longrightarrow}E\left(p\left(u_0-B_s\right)q\left(v_0-B_t\right)\right) \]
for almost all $\left(s,t\right)\in\left[0,1\right]^2$. This is fact holds true for 
\[ \left(s,t\right)\notin S_0:=\left\{s=0\right\}\cup\left\{t=0\right\}\cup\left\{s=t\right\}. \] 
Indeed, we have
\begin{equation}\label{eq:expectation_of_shifted_functions}
E\left(p\left(u-B_s\right)q\left(v-B_t\right)\right)=\iint_{\reals^{2d}}p\left(x\right)q\left(y\right)f_{s,t}\left(x-u,y-v\right)\ud x\ud y.
\end{equation}
The function $f_{s,t}\left(x-u,y-v\right)$ is continuous in $\left(x,y,u,v\right)$ for each fixed $\left(s,t\right)\notin S_0$. In addition, if $\abs{u},\abs{v}\leq M$, then for some suitable constant $\tilde{M}$ (which depends only on $M$),
\[ f_{s,t}\left(x-u,y-v\right)\leq f_{s,t}\left(x,y\right)e^{\tilde{M}\norm{\Sigma_{s,t}^{-1}}\left(1+\abs{x}+\abs{y}\right)} \]
($\Sigma_{s,t}$ here denotes the covariance matrix of $\left(B_s,B_t\right)$). This provides dominated convergence in~\eqref{eq:expectation_of_shifted_functions} since then for all $\left(s,t\right)\notin S_0$,
\[ \abs{p\left(x\right)}\abs{q\left(y\right)}f_{s,t}\left(x-u,y-v\right)\leq C^2f_{s,t}\left(x,y\right)e^{\alpha\left(\abs{x}+\abs{y}\right)+\tilde{M}\norm{\Sigma_{s,t}^{-1}}\left(1+\abs{x}+\abs{y}\right)} \]
which is in $L^1\left(\reals^{2d}\right)$.
\end{proof}

\bigskip

We now proceed to the proof of the Theorem.
\begin{proof}[Proof of Theorem~\ref{theorem:main}]
Note that since $F$ and $G$ are locally bounded by assumption~\eqref{eq:exponential_growth_assumption}, the processes $u^t_r:=F\left(B_{r}\right)\One{0}{t}$ and $v^s_r:=G\left( B_{r}\right)\One{0}{s}$ belong a.s. to $\abs{\rkhs}^d$.

We will work our way from regular $F,G$'s to the general case.\medskip\\
1. $F,G\!\in\!C^{\infty}_{\operatorname{c}}(\mathbb{R}^d;\mathbb{R}^d)$.\medskip\\
Fix $s,t\in\left[0,T\right]$. The following Lemma identifies the random Hilbert-Schmidt operators appearing in~\eqref{absfor:reg}.
\begin{lemma}\label{prop:derivative}
  Assume that $\varphi\!\in\!C^2\left(\mathbb{R}^d\right)$ has bounded first and second derivatives, and for some $r\!\in\![0,T]$ denote $u_s=\varphi\left(B_s\right)\mathbf{1}_{\left[0,r\right]}(s),\,\ 0\leq s\leq T$. Then $u\!\in\!\mathbb{D}^{1,2}\left(\mathcal{H}\right)$ (in particular $u\!\in\!\operatorname{Dom}{\delta}$)\ and $D^iu$ is represented by the following kernel in the sense of~\eqref{kernel_representation}:
\begin{equation}\label{grad_ker_rep}
  k\left(s,t\right)=\partial_i\varphi\left(B_s\right)\One{t}{r}\left(s\right),\ \ t\leq r\ \text{(and zero otherwise)}.
\end{equation}
\end{lemma}
Note that the assertion of the Lemma can be written explicitly as (assuming $t\leq r$)
\begin{equation*}
    (D^iu)h(t)=\alpha_H\int_0^T\ud \tau\int_t^rh\left(\tau\right)\abs{\tau-\sigma}^{2H-2}
              \partial_i\varphi\left(B_{\sigma}\right)\ud \sigma,\ \ h\in\abs{\mathcal{H}}.
\end{equation*}
By Lemma~\ref{prop:derivative}\ \ (with $\varphi\!=\!F,G$ and $r\!=\!t,s$ respectively)\,\
$D^i\left(u_j\mathbf{1}_{\left[0,t\right]}\right)$ and $D^j\left(v_i\mathbf{1}_{\left[0,s\right]}\right)$ are represented in the sense of~\eqref{kernel_representation} by the kernel~\eqref{grad_ker_rep} . We use these kernels for the derivatives in~\eqref{multi_dim_abs_for}, and Proposition~\ref{lemma:kernels} to evaluate the trace:
 \begin{align*}
   E\left(\int_0^tF\left( B_{\tau}\right)\cdot\ud B_{\tau}\,\int_0^sG\left( B_{\sigma}\right)\cdot\ud
                                                                             B_{\sigma}\right)&=\sum_{i=1}^dE\left<F_i\left(B_{\cdot}\right)\One{0}{t},G_i\left(B_{\cdot}\right)\One{0}{s}\right>_{\mathcal{H}}+\\
        &\hspace{-5cm}+\alpha_H^2\int_0^t\ud{\tau_1}\int_0^s\ud{\sigma_1}\int_{{\tau_1}}^{t}\ud{\tau_2}\int_{{\sigma_1}}^s
                \sum_{i,j=1}^dE\left(\partial_iF_j\left(B_{{\tau_2}}\right)\partial_jG_i\left(B_{{\sigma_2}}\right)\right)\cdot\\
                &\hspace{0.5cm}\cdot\abs{{\tau_2}-{\sigma_1}}^{2H-2}\abs{{\tau_1}-{\sigma_2}}^{2H-2}\ud{\sigma_2}.
 \end{align*}
The expression in the Theorem now follows by changing the order of integration (by where the integrals according to $\tau_1$ and $\sigma_1$ are carried out first, resulting in $\gamma(\tau_2,\sigma_2)$).\medskip\\
2. $F,G\!\in\!C^{\infty}(\mathbb{R}^d;\mathbb{R}^d)$ satisfy Assumption~\eqref{exponential_growth} and in addition \\
$\int_{\reals^{2d}}L\left(x,y\right)\abs{DF\left(y\right)}\abs{DF\left(x\right)}\ud x\!\ud y<\infty$ (and the same for $G$). 
\medskip\\
We begin with $F$. Let $\psi_n$ be a $C^{\infty}$ bump function on $\abs{x}\leq n$: $\psi_n=1$ on $\abs{x}\leq n$, $\psi_n=0$ outside of $\abs{x}\leq n+1$ and $0\leq\psi_n\leq 1$ in-between. Set $F_n=F\psi_n$. Define $u_n^{t}=F_n\left(B_{\cdot}\right)\mathbf{1}_{\left[0,t\right]}$. To see that $\left(u_n^t\right)_{n=1}^{\infty}$ converges to $u^t$ in $L^2\left(\mathbb{P};\mathcal{H}^d\right)$, note that $u_n^t\to u^t$ almost everywhere as $n\to\infty$. Now by Assumption~\eqref{exponential_growth},
 \begin{equation*}
   \abs{\left(u_n^t\left({\tau}\right)-u^t\left({\tau}\right)\right)\left(u_n^t\left({\sigma}\right)-u^t\left({\sigma}\right)\right)\abs{{\sigma}-{\tau}}^{2H-2}}\!\leq\!
   C^2e^{\beta\left(\abs{B_{\tau}}+\abs{B_{\sigma}}\right)}\abs{{\sigma}-{\tau}}^{2H-2}\ .
  \end{equation*}
 The right-hand-side is  integrable in $\left[0,T\right]^2$; thus, by the dominated convergence theorem, $\lim_{n\rightarrow\infty}\norm{u_n^t-u^t}^2_{\mathcal{H}}=0$. Since
  \begin{equation*}
   \norm{u^t_n-u^t}_{\mathcal{H}}^2\leq\alpha_HM^2\int_{\left[0,T\right]^2}
C^2e^{\beta\left(\abs{B_{\tau}}+\abs{B_{\sigma}}\right)}\abs{\tau-\sigma}^{2H-2}\ud\tau\ud\sigma\in L^1\left(\mathbb{P}\right),
  \end{equation*}
  it follows (from dominated convergence) that $\lim_{n\rightarrow\infty}E\norm{u_n^t-u^t}^2_{\mathcal{H}}=0$.

  The next step is to show that $\left(\delta u_n^t\right)_{n=1}^{\infty}$ is a Cauchy sequence in $L^2\left(\mathbb{P}\right)$. By Remark~\ref{multi_dim_malliavin},
  \begin{equation*}
    E\abs{\delta u_n^t-\delta u_m^t}\leq\sum_{i=1}^dE\abs{\delta^iu_n^{i,t}-\delta^iu_m^{i,t}}.
  \end{equation*}
We can now deduce from Formula~\eqref{form}, seeing as the $F_n$'s are $C^{\infty}_{\operatorname{c}}$, that
  \begin{align*}
   E\abs{\delta^i u_n^{i,t}-\delta^i u_m^{i,t}}^2=&E\norm{u_n^{i,t}-u_m^{i,t}}^2_{\mathcal{H}}+\\
   &\hspace{-3cm}+\alpha_H\int_{\left[0,t\right]^2}\!\!\!\!\gamma\left(\tau,\sigma\right) E\left[\left(\partial_iF^i_n\left(B_{\tau}\right)-\partial_iF^i_m\left(B_{\tau}\right)\right)\left(\partial_iF^i_n\left(B_{\sigma}\right)-\partial_iF^i_m\left(B_{\sigma}\right)\right)\right]\ud\tau\ud\sigma.
  \end{align*}
  By Fubini (and recalling $L\left(x,y\right)$'s definition in~\eqref{L}), it remains to prove that
  \begin{equation}\label{Cauchy}
   \lim_{n,m\to\infty}\int_{\mathbb{R}^{2d}}\left(\partial_iF^i_n\left(x\right)-\partial_iF^i_m\left(x\right)\right)\left(\partial_iF^i_n\left(y\right)-\partial_iF^i_m\left(y\right)\right)L\left(x,y\right)\ud x\ud y=0.
  \end{equation}
  (To be precise, this is required for $t$ instead of $T$ in~\eqref{L}, but this is of course inconsequent.) This would follow from
    \begin{equation}\label{Cauchy2}
   \lim_{n,m\to\infty}\!\int_{\mathbb{R}^{2d}}\!\!\partial_iF^i_n\left(x\right)\partial_iF^i_m\left(y\right)L\left(x,y\right)\ud x\ud y=\int_{\mathbb{R}^{2d}}\!\!\partial_iF^i\left(x\right)\partial_iF^i\left(y\right)L\left(x,y\right)\ud x\ud y.
  \end{equation}
  Since $\partial_iF^i_m=\partial_i\psi_mF^i+\psi_m\partial_iF^i$, there are four terms on the left-hand-side of~\eqref{Cauchy2}. Three of them tend to zero, and the fourth to the right-hand-side. Indeed, these all follow easily from dominated convergence considering the assumptions on $F$.

  Therefore $\left(\delta u_n^t\right)_{n=1}^{\infty}$ is a Cauchy sequence in $L^2\left(\mathbb{P}\right)$. Since $\delta$ is a closed operator, it follows that $u^t\in\operatorname{Dom}\left(\delta\right)$ and that
 \begin{equation*}
  \delta u^t=\lim_{n\rightarrow\infty}\delta u_n^t\qquad\text{in
  }L^2\left(\mathbb{P}\right).
 \end{equation*}
All of the above holds for $v_n^s=\left(G\psi_n\right)\left(B_{\cdot}\right)\One{0}{s}$ as well, and
 \begin{align*}
  E\delta u^t\delta v^s=\lim_{n\rightarrow\infty}E\delta u_n^t\delta v_n^s&\stackrel{\text{\eqref{form}}}{=}\int_0^t\int_0^sM\left(\tau,\sigma\right)\abs{\tau-\sigma}^{2H-2}\ud\tau\ud\sigma\\
  &\hspace{-0.6cm}+\alpha_H\lim_{n\rightarrow\infty}\int_{\mathbb{R}^{2d}}\sum_{i,j=1}^d\partial_iF^j_n\left(x\right)\partial_jG^i_n\left(y\right)L^{t,s}\left(x,y\right)\ud
  x\ud y,
 \end{align*}
 where
 \begin{equation*}
  L^{t,s}\left(x,y\right)=\int_0^t\int_0^s\gamma\left(\tau,\sigma\right)f_{\tau,\sigma}\left(x,y\right)\ud\tau\ud\sigma.
 \end{equation*}
 Finally, just like Equation~\eqref{Cauchy2}, we have
  \begin{equation}\label{replacable_2}
  \lim_{n\rightarrow\infty}\!\iint_{\mathbb{R}^{2d}}\!\!\partial_iF^j_n\left(x\right)\partial_jG^i_n\left(y\right)L^{t,s}\left(x,y\right)\ud  x\ud y=\int_{\mathbb{R}^{2d}}\!\!L^{t,s}\left(x,y\right)\partial_iF^j\left(\ud x\right)\partial_jG^i\left(\ud y\right).
 \end{equation}

3. $F,G$ satisfying the assumptions of the Theorem.\medskip \\
All that is now required is the following Lemma:
\begin{lemma}\label{convolution_approximation}
	Let $\varphi_n:\reals^d\to\reals$ be an approximation of the identity: $\varphi_n\left(x\right)=n\varphi\left(nx\right)$ where $\varphi$ is a non-negative $C^{\infty}_{\operatorname{c}}$ function, supported in the closed unit ball, such that $\int_{\reals^d}\varphi\left(x\right)\ud x=1$. Set $F_n=F\ast\varphi_n$ (similarly for $G_n$). Then
	 \begin{equation}\label{conv_limit}
  \lim_{m,n\rightarrow\infty}\iint_{\mathbb{R}^{2d}}\partial_iF^j_m\left(x\right)\partial_jG^i_n\left(y\right)L\left(x,y\right)\ud  x\ud y=\int_{\mathbb{R}^{2d}}L\left(x,y\right)\partial_iF^j\otimes\partial_jG^i\left(\ud x,\ud y\right).
 \end{equation}
\end{lemma}
As in the previous step, set $u_n^{t}=F_n\left(B_{\cdot}\right)\mathbf{1}_{\left[0,t\right]}$, and basically repeat the same pattern; first, note that $F_n$ satisfies the assumptions in the previous step by Lemma~\ref{convolution_approximation} and since the $F_n$'s inherit from $F$ the same exponential growth assumption~\eqref{exponential_growth}, with (possibly different) constants $C$ and $\beta$ (which do not depend on $n$). Now, $u_n^t\to u^t$ almost everywhere as $n\to\infty$ - this follows for example from~\cite[Theorem~8.15]{Folland}. The remaining equations hold true for the same reasons, with the exception of Equations~\eqref{Cauchy2} and~\eqref{replacable_2} which now follow from Lemma~\ref{convolution_approximation}.

This concludes the proof of Theorem~\ref{theorem:main}.
\end{proof}
\bigskip

We end this section with the proofs of Lemmas~\ref{prop:derivative} and~\ref{convolution_approximation}.
\begin{proof}[Proof of Lemma~\ref{prop:derivative}]		
To simplify the proof, essentially without affecting it otherwise, we let $r\!=\!T$.

\noindent
Recall that $u\in\left|\mathcal{H}\right|$ almost surely. We shall approximate $u$ (in $L^2\left(\mathbb{P},\mathcal{H}\right)$) by smooth random elements in such a way that the resulting sequence of derivatives converges (in $L^2\left(\mathbb{P},\mathcal{H}\otimes\mathcal{H}\right)$). For each $n\in\mathbb{N}$, denote $t_k=t_k^n=\frac{kT}{2^n}$, $k=0,\ldots,2^n$ and
  \[ u^n=\sum_{k=1}^{2^n}\varphi\left(B_{t_{k-1}}\right)\mathbf{1}_{\left[t_{k-1},t_k\right]}\ \ \
                                                               \in L^{2}\left(\mathbb{P};\mathcal{H}\right)\,. \]
 By definition\ \ \ \ \
 $\displaystyle{D^iu^n\!=\!\sum_{k=1}^{2^n}\partial_i\varphi\!\left(B_{t_{k-1}}\right)
        \mathbf{1}_{\left[0,t_{k-1}\right]}\otimes\mathbf{1}_{\left[t_{k-1},t_k\right]}}$\,,\ \ \ \
 so that for any $h\in\abs{\mathcal{H}}$
\begin{equation}  \label{Dun}
D^iu^n\left(h\right)\left(x\right)
   =\alpha_H\sum_{k=1}^{2^{n}}\int_0^TH\left(t\right)
              \mathbf{1}_{\left[t_{k-1},t_{k}\right]}\left(t\right)
              \partial_i\varphi\!\left(B_{t_{k-1}}\right)\mathbf{1}_{\left[x,T\right]}\left(t_{k-1}\right)\ud t,
\end{equation}
with\ \ \
$\displaystyle{H\left(t\right)\!=\!\!\int_0^T\!\!h\left(s\right)\abs{s\!-\!t}^{2H-2}\!\ud s}$\,.\ \
The lemma will thus follow once we show that $u_n\!\to\!u$\ in\ $L^{2}\left(\mathbb{P},\mathcal{H}\right)$ and that $\left\{ D^iu^{n}\right\} _{n=1}^{\infty}$ is a Cauchy sequence in
$L^{2}\left(\mathbb{P};\mathcal{H}\!\otimes\!\mathcal{H}\right)$; indeed, this implies that $u\in\mathbb{D}^{1,2}\left(\mathcal{H}\right)$\,,\ $u^n\to u$ in $\mathbb{D}^{1,2}$\ and in particular\ \ $Du^n\to Du$. The equality~\eqref{grad_ker_rep} will then follow directly from~\eqref{Dun}.\medskip

 \noindent
 Concerning the first assertion
\begin{equation*}
u_{t}-u_{t}^{n}=\sum_{k=1}^{2^{n}}\left[\varphi\left(B_{t}\right)-\varphi\left(B_{t_{k-1}}\right)\right]\mathbf{1}_{\left[t_{k-1},t_{k}\right]}
          \ \ :=\ \ \sum_{k=1}^{2^{n}}v_t^{n,k};
\end{equation*}
whence, by the H\"{o}lder inequality,
\begin{equation*}
 E\left\Vert u-u^{n}\right\Vert_{\mathcal{H}}^{2}\
        \leq\ 2^{n}\sum_{k=1}^{2^n}E\left\Vert v^{n,k}\right\Vert_{\mathcal{H}}^{2}\
        \leq\ 2^{2n}d\norm{\nabla \varphi}_{\infty}^{2}\left|\frac{T}{2^{n}}\right|^{4H}\
        \goes\ 0.
\end{equation*}
As for the sequence of derivatives, denote $s_{i}=t_i^{n+1}=\frac{Ti}{2^{n+1}}$, \ \ and then
\begin{equation*}
u^{n}-u^{n+1}=\sum_{i=1}^{2^{n}}\left[\varphi\left(B_{s_{_{2i-2}}}\right)-\varphi\left(B_{s_{_{2i-1}}}\right)\right]
                                                             \mathbf{1}_{\left[s_{_{2i-1}},s_{_{2i}}\right]},
\end{equation*}
and
\begin{eqnarray*}
 \lefteqn{\left\Vert D^iu^{n}-D^iu^{n+1}\right\Vert ^{2}\leq}\\
   &&\leq 2^{n}\left(\frac{T}{2^{n+1}}\right)^{2H}
   \sum_{i=1}^{2^{n}}\left\Vert\partial_i\varphi\left(B_{s_{_{2i-2}}}\right)\mathbf{1}_{\left[0,s_{_{2i-2}}\right]}
                 -\partial_i\varphi\left(B_{s_{_{2i-1}}}\right)\mathbf{1}_{\left[0,s_{_{2i-1}}\right]}\right\Vert^{2}.
\end{eqnarray*}
Since $\left[0,s_{_{2i-1}}\right]=\left[0,s_{_{2i-2}}\right]\cup\left[s_{_{2i-2}},s_{_{2i-1}}\right]$ is a non-overlapping union,
\begin{eqnarray*}
  \lefteqn{\left\Vert\partial_i\varphi\left(B_{s_{_{2i-2}}}\right)\mathbf{1}_{\left[0,s_{_{2i-2}}\right]}
      -\partial_i\varphi\left(B_{s_{_{2i-1}}}\right)\mathbf{1}_{\left[0,s_{_{2i-1}}\right]}\right\Vert^{2}}\\
  &&\leq 2\Bigg(\left[\partial_i\varphi\!\left(\!B_{s_{_{2i-2}}}\!\right)
             -\partial_i\varphi\!\left(\!B_{s_{_{2i-1}}}\!\right)\right]^{2}\left|s_{_{2i-2}}\right|^{2H}
                +\partial_i\varphi\left(B_{s_{_{2i-1}}}\right)^{2}\left(\frac{T}{2^{n+1}}\right)^{2H}\Bigg)\\
  &&\leq T^{2H}\left(2\|\nabla\partial_i\varphi\|_\infty^2\left(B_{s_{_{2i-2}}}\!\!\!\!\!-\!B_{s_{_{2i-1}}}\right)^2
                      +\|\partial_i\varphi\|_\infty^2\,2^{-2Hn}\right).
\end{eqnarray*}
  It follows that $E\left\Vert Du^{n}-Du^{n+1}\right\Vert ^{2}\leq C_1\,2^{(2-4H)n}$ for some constant $C_1$ which depends on $H,d,T$ and $\varphi$ but not on $n$. Thus for any $n_0$ and all $n\!>\!m\!\ge\!n_0$
 \begin{align*}
   \left\Vert D^iu^{n}-D^iu^{m}\right\Vert_{L^{2}\left(\mathbb{P};\mathcal{H}\otimes\mathcal{H}\right)}
      &\leq \sum_{k=m}^{n-1}\left\Vert D^iu^k-D^iu^{k+1}\right\Vert_{L^{2}\left(\mathbb{P};\mathcal{H}\otimes\mathcal{H}\right)}\\
      &\leq C_2\,\sum_{k=m}^{n-1}2^{\left(1-2H\right)k}
      \leq C_3\, 2^{\left(1-2H\right)n_0}
  \end{align*}
  (for suitable $C_2,C_3$) so that $\left\{ D^iu^{n}\right\}_{n=1}^{\infty}$ is a Cauchy sequence in $L^{2}\left(\mathbb{P},\mathcal{H}\!\otimes\!\mathcal{H}\right)$.
\end{proof}

\begin{proof}[Proof of Lemma~\ref{convolution_approximation}]
The fact that $\partial_iF^j_m\left(x\right)\partial_jG^i_n\left(y\right)L\left(x,y\right)\in L^1\left(\reals^{2d}\right)$ follows in essentially the same way as the calculation below.

Recall that $L\ast\left(\partial_iF^j\otimes\partial_jG^i\right)$ is continuous at the origin. Now compute 
\begin{align*}
\int_{\mathbb{R}^{2d}}\partial_iF^j_m&\left(x\right)\partial_jG^i_n\left(y\right)L\left(x,y\right)\ud  x\ud y\\
&=\int_{\mathbb{R}^{2d}}\partial_iF^j\ast\varphi_m\left(x\right)\partial_jG^i\ast\varphi_n\left(y\right)L\left(x,y\right)\ud  x\ud y\\
&=\int_{\mathbb{R}^{4d}}\varphi_m\left(x-x'\right)\varphi_n\left(y-y'\right)L\left(x,y\right)\partial_iF^j\otimes\partial_jG^i\left(\ud x',\ud y'\right)\ud  x\ud y\\
&=\int_{\mathbb{R}^{4d}}\varphi_m\left(\tilde{x}\right)\varphi_n\left(\tilde{y}\right)L\left(\tilde{x}+x',\tilde{y}+y'\right)\partial_iF^j\otimes\partial_jG^i\left(\ud x',\ud y'\right)\ud  \tilde{x}\ud \tilde{y}\\
&=\int_{\reals^{2d}}L\ast\left(\partial_iF^j\otimes\partial_jG^i\right)\left(-x,-y\right)\varphi_m\left(x\right)\varphi_n\left(y\right)\ud x\ud y,
\end{align*}
since $L\left(-x,-y\right)=L\left(x,y\right)$.
Assume without loss of generality that $m\geq n$. Then in the last integral above, $\abs{x},\abs{y}\leq \frac{1}{n}$; by choosing $n$ large enough such that 
\begin{equation*}
\abs{L\ast\left(\partial_iF^j\otimes\partial_jG^i\right)\left(-x,-y\right)-L\ast\left(\partial_iF^j\otimes\partial_jG^i\right)\left(0,0\right)}<\epsilon,
\end{equation*}
we now conclude the proof, since the second term above is equal to the right-hand-side of~\eqref{conv_limit} and by the properties of the $\varphi$'s.
\end{proof}
\noindent


\section{An Application to the Fractional Bessel Process}   \label{Bessel}
When $d\!>\!1$ the $d$-dimensional fractional Bessel process\
 $R_t=\left|B_t\right|$\ \
satisfies
\begin{equation}\label{bessel_multi}
  R_t=\int_0^t\frac{B_s}{R_s}\cdot\ud B_s+H\left(d-1\right)\int_0^t\frac{s^{2H-1}}{R_s}\ud s.
\end{equation}
This follows from the It\^{o} formula for fractional Brownian motion (see for example~\cite{Guerra&Nualart}).
For $d\!=\!1$ one has an analogue of Tanaka's formula (see~\cite{Hu&Oksendal&Salopek}):
\begin{equation}\label{tanaka}
  R_t=\int_0^t\sgn\left(B_s\right)\ud B_s+\text{local time process}.
\end{equation}
In~\cite{Hu_Nualart} Y. Hu and D. Nualart asked whether, by analogy with Brownian motion, the $d$-dimensional stochastic integral  $X_t\!=\!\!\displaystyle{\int_0^t}\frac{B_s}{R_s}\cdot\ud B_s$\ appearing in~\eqref{bessel_multi} and~\eqref{tanaka} is itself a $d$-dimensional fractional Brownian motion, 
 and subsequently provided a negative answer using Wiener chaos expansions.
We now show (when $H\!>\!\nicefrac{1}{2}$) how this conclusion can be arrived at directly from the covariance function of $X_t$ obtained in Theorem~\ref{generalization}.

\noindent
Let then $F\left(x\right)\!=\!G\left(x\right)\!=\!\hat{x}:=\bfrac{x}{\abs{x}}$ in Theorem~\ref{generalization}, whose conditions indeed hold for $d=1$, as already observed in the discussion following its statement, in particular~(\ref{finite.integral}). The case $d\geq 2$ follows shortly. For $d=1$, and since $F'\!=\!G'\!=\!2\delta_0$, the terms in~\eqref{form:general} are given by
  \begin{eqnarray}\label{eq_sgn.der}
    4\int_{\reals^2}L_{t,s}\left(x,y\right)\delta_0\left(\ud x\right)\delta_0\left(\ud y\right)&=&\int_0^t\int_0^s\frac{\gamma\left(\tau,\sigma\right)}{\pi\tau^H\sigma^H\sqrt{1-\rho^2\left(\tau,\sigma\right)}}\ud\tau\ud\sigma,\\ \label{Msgn}
      \E\left(\sgn\left(B_{\tau}\right)\sgn\left(B_{\sigma}\right)\right)&=&\frac{2}{\pi}\arccos\left(\sqrt{1-\rho^2\left(\tau,\sigma\right)}\right).
  \end{eqnarray}
  where
  \begin{equation}\label{corell}
   \rho\left(\tau,\sigma\right)=\frac{R_H\left(\tau,\sigma\right)}{\tau^H\sigma^H}\ .\\
     \end{equation}
  (The standard formula for the first quadrant Gaussian measure was used to obtain~\eqref{Msgn}.) It doesn't seem possible to explicitly compute the resulting expression for $\E\left(\int_0^t\sgn B_\tau\ud B_\tau\int_0^s\sgn B_\sigma\ud B_\sigma\right)$ (recall $\gamma(t,s)'s$ definition in~\eqref{gamma}):
  \begin{equation}\label{sgncovar}
         \frac{1}{\pi}\!\!\!\int_0^t\!\!\!\int_0^s\!\left(\alpha_H\frac{2}{\pi}\arccos\left(\sqrt{1-\rho^2\left(\tau,\sigma\right)}\right)\abs{\tau\!-\!\sigma}^{2H-2}\!
                                              +\frac{\gamma\left(\tau,\sigma\right)}{\tau^H\sigma^H\sqrt{1-\rho^2\left(\tau,\sigma\right)}}\right)\!\ud\sigma\!\ud\tau,
  \end{equation}
   however it is easy to see by substitution that its mixed second derivative does not coincide with
  \begin{eqnarray*}
    \frac{\partial^2R_H(t,s)}{\partial t\partial s}=\alpha_{H}\abs{t-s}^{2H-2}.
  \end{eqnarray*}
  \noindent
  Thus, unlike the case $H\!\!=\!\!\nicefrac{1}{2}$,\ \
  $X_t\!\!=\!\!\int_0^t\sgn B_\tau\ud B_\tau$ is not itself a fractional Brownian motion with parameter $H$ or in fact with any other parameter $\widetilde{H}$.

For $d\geq 2$, assuming for the moment that Condition~\ref{condition:continuity_of_convolution} in Theorem~\ref{theorem:main} holds true (see below), we
 obtain:
\begin{eqnarray}
  \E X_tX_s\ &\!\!\!\!\!-\!\!\!\!\!&\alpha_H\!\!\!\int_0^t\!\!\!
                    \int_0^s\abs{\tau-\sigma}^{2H-2}\,E\widehat{B_{\tau}}\!\cdot\!\widehat{ B_{\sigma}}\,d\sigma d\tau\nonumber\\
                &&\hspace{-1cm}\!\!\!\!\!=\alpha_H\,\gamma(\tau,\sigma)
                   \int_0^t\!\!\!\int_0^s\sum_{i,j=1}^d
                        E\!\left(\frac{\norm{B_{\tau}}^2\delta_{ij}-B^i_{\tau}B^j_{\tau}}{\norm{B_{\tau}}^3}\right)
                   \left(\frac{\norm{B_{\sigma}}^2\delta_{ij}-B^i_{\sigma}B^j_{\sigma}}{\norm{B_{\sigma}}^3}\right)
                                                      \!\ud\sigma\!\ud\tau\nonumber\\
                &\!\!\!\!\!=\!\!\!\!\!\!&\alpha_H\,\gamma(\tau,\sigma)
                                                            \int_0^t\!\!\!\int_0^s
                E\!\left(\frac{(d\!-\!2)+(\widehat{B_\tau}\cdot\widehat{B_\sigma})^2}
                             {\norm{B_\tau}\norm{B_\sigma}}\right)\ud\sigma\!\ud\tau
       \label{extxs}
\end{eqnarray}
If $X$ were a fBm with parameter $H$ (by self-similarity arguments it cannot be a fBM with any other parameter), by taking the mixed second derivative of~(\ref{extxs}) and dividing both its sides by $\alpha_H|t\!-\!s|^{2H-2}$ we obtain	 
 \begin{equation}
  1-E\widehat{B_t}\cdot\widehat{B_s}
     =\frac{\gamma(\tau,\sigma)}{|t\!-\!s|^{2H-2}}\,
            E\!\left(\frac{(d\!-\!2)+(\widehat{B_t}\cdot\widehat{B_s})^2}{\norm{B_t}\norm{B_s}}\right)
     \le \frac{(d\!-\!1)\gamma(\tau,\sigma)}{|t\!-\!s|^{2H-2}}\,
            E\!\left(\frac{1}{\norm{B_t}\norm{B_s}}\right)
 \end{equation}
 Fix $s\!>\!0$. From~(\ref{gamma})\,, $\gamma(t,s)=O(t^{4H-3})$ as $t\!\to\!\infty$. We moreover claim that
 \begin{equation} \label{asymptotics}  		
   E\!\left(\frac{1}{\norm{B_t}\norm{B_s}}\right)=O(t^{-H})\hspace{1cm}{\rm as}\ t\!\to\!\infty
 \end{equation}
 from which it would follow that $\displaystyle{\lim_{t\to\infty}}E\widehat{B_t}\!\cdot\!\widehat{B_s}=1$  and thus $\displaystyle{E\widehat{B_{s}}\!\cdot\!\widehat{B_{s'}}=1}$ for any $s,s'$ (as can be seen from the estimate\ \ $1\!-\!\hat{x}\!\cdot\!\hat{y}
     \le 2\left(\rule[-4pt]{0pt}{10pt}(1\!-\!\hat{x}\!\cdot\!\hat{z})+(1\!-\!\!\hat{y}\!\cdot\!\hat{z})\right)$
 for all\ $x,y,z \in\!\Reals^d$)\,. However $(B_s,B_{s'})$ is clearly non-degenerate for $s\!\ne\!s'$\ and
 this contradiction will show, once~(\ref{asymptotics}) is verified, that $X_t$ cannot be a fBM. \medskip

 \noindent
 For $d\!\ge\!3$, (\ref{asymptotics}) follows just by the Cauchy Schwartz inequality and self-similarity:
  \[ E\!\left(\frac{1}{\norm{B_t}\norm{B_s}}\right)\le
     t^{-H}E\!\left(\!\frac{1}{\norm{B_1}^2}\!\right)^{\!\!1/2}\!\!E\!\left(\!\frac{1}{\norm{B_s}^2}\!\right)^{\!\!1/2}\]
 whereas for $d\!=\!2$\ (in which case $E\bfrac{1}{\|B_s\|^2}\!=\!\infty$) 
 we denote $(B_t,B_s)\sim {\rm N}(0,\Sigma^{(2)})$ where
 $\Sigma^{(2)}\!\!=\!\!\left(\!\!\!
               \begin{array}{cc}
                 t^{2H}I_2 & R_H(t,s)I_2 \\
                 R_H(t,s)I_2 & s^{2H}I_2 \\
               \end{array}\!\!\!\right)$\,,\ namely the Kroeneker product $\Sigma\otimes I_2$ of the scalar covariance matrix
 and the $2\mbox{\scriptsize $\times$}2$ identity matrix. Note that $|\Sigma^{(2)}|\!=\!|\Sigma|^2$ (here $|\cdot|=\det$)
 and that $\left(\Sigma^{(2)}\right)^{-1}\!=\!\left(\Sigma^{-1}\right)^{(2)}$. Then, passing to polar coordinates $(r,\theta,r',\theta')$,

 \begin{align}
  E\,\frac{1}{\|B_t\|\,\|B_s\|}
        &=\frac{1}{4\pi^2|\Sigma|}\int_{[0,2\pi]^2}\hspace{-16pt}d\theta d\theta'
                                     \int_{[0,\infty)^2}\hspace{-16pt}\exp\frac{\Sigma^{-1}_{11}r^2
                                               \!+\!2\Sigma^{-1}_{12}\cos(\theta\!-\!\theta')rr'\!+\!\Sigma^{-1}_{22}r'^2}{2}\ drdr'\nonumber\\
        &=\rule{0pt}{25pt}\frac{1}{4\pi^2|\Sigma|}\int_{[0,2\pi]^2}\hspace{-6pt}
          2\pi|\Sigma_{\theta,\theta'}|^{1/2}\,P_{\theta,\theta'}(R\!\ge\!0,\,R'\!\ge\!0)\,d\theta d\theta'\nonumber\\
          &\le\rule{0pt}{25pt}\frac{2\pi}{|\Sigma|^{1/2}} \label{minusonemoment}
 \end{align}
 \smallskip

 \noindent
     where $P_{\theta,\theta'}\sim{\rm N}(0,\Sigma_{\theta,\theta'})$\ \ with\ \
      $\left(\Sigma_{\theta,\theta'}\right)^{-1}\!=\!\left(\!\!\begin{array}{ll}
       \Sigma^{-1}_{11} & \Sigma^{-1}_{12}\cos(\theta\!-\!\theta') \\
       \Sigma^{-1}_{21}\cos(\theta\!-\!\theta')& \Sigma^{-1}_{22}\end{array}\!\!\right)$\ \ \ and the last inequality follows from
       $\left|\left(\Sigma_{\theta,\theta'}\right)^{-1}\right|\!\ge\!\left|\Sigma^{-1}\right|$.
    The obvious asymptotic estimate\ \ $|\Sigma|=\left(s^{2H}\!+\!o(1)\right)\,t^{2H}$\ as $t\!\to\!\infty$\,, applied to~(\ref{minusonemoment}),\ proves~(\ref{asymptotics}) for $d\!=\!2$.
    
    It remains to be seen that Condition~\eqref{condition:continuity_of_convolution} in Theorem~\ref{theorem:main} indeed holds true for $d\geq 2$. For $d\geq 3$ we denote $p_{t,s}\left(u,v\right)=\E \partial_iF^j\left(B_s-u\right)\partial_kF^l\left(B_t-v\right)$. Then
\begin{equation}\label{example_condition_1}
    	L\ast\partial_iF^j\otimes\partial_kF^l\left(u,v\right)=\int_0^{\tau}\int_0^{\sigma}p_{t,s}\left(u,v\right)\gamma\left(t,s\right)\ud t\ud s.
\end{equation}
Evaluating:	
    \begin{align*}
    	\abs{p_{t,s}\left(u,v\right)}&\leq\E\left(\frac{1}{\abs{B_s-u}}\frac{1}{\abs{B_t-v}}\right)\leq \sqrt{\E\left(\frac{1}{\abs{B_s-u}^2}\right)\E\left(\frac{1}{\abs{B_t-v}^2}\right)}\\
	\text{\scriptsize\fbox{Anderson's Lemma}} & \leq \sqrt{\E\left(\frac{1}{\abs{B_s}^2}\right)\E\left(\frac{1}{\abs{B_t}^2}\right)}\stackrel{\text{\scriptsize\fbox{Self-similarity}}}{=}\frac{C}{s^Ht^H},
    \end{align*}
it follows that we have dominated convergence in~\eqref{example_condition_1}. For $d=2$, dominated convergence in~\eqref{example_condition_1} already follows from~\eqref{minusonemoment} and~\eqref{finite.integral}. Finally, it remains to check for the continuity of $p_{t,s}$ at the origin (for almost all $s,t$). We can write, denoting $h=\partial_iF^j\otimes\partial_kF^l$,
\[ p_{t,s}\left(u,v\right)=C_{t,s}\int_{\reals^{2d}}h\left(x-u,y-v\right)e^{-\frac{1}{2}\left(A_{t,s}\abs{x}^2+B_{t,s}x\cdot y+D_{t,s}\abs{y}^2\right)}\ud x\ud y. \]
By a simple change of variables, we now have
\begin{align}\label{example_condition_2}
	\abs{p_{t,s}\left(u,v\right)-p_{t,s}\left(0,0\right)}&\\
	&\hspace{-3cm}\leq\int_{\reals^{2d}}h\left(x,y\right)f_{t,s}\left(x,y\right)\abs{e^{-\frac{1}{2}\left(A_{t,s}\left(2x\cdot u+\abs{u}^2\right)+B_{t,s}\left(x\cdot v+u\cdot y+u\cdot v\right)+D_{t,s}\left(2y\cdot v+\abs{v}^2\right)\right)}-1}\ud x\ud y\nonumber
\end{align}
For $\abs{u},\abs{v}\leq 1$,
\[ \abs{e^{-\frac{1}{2}\left(A_{t,s}\left(2x\cdot u+\abs{u}^2\right)+B_{t,s}\left(x\cdot v+u\cdot y+u\cdot v\right)+D_{t,s}\left(2y\cdot v+\abs{v}^2\right)\right)}-1}\leq e^{A_{t,s}\abs{x}+B_{t,s}\left(\abs{x}+\abs{y}+1\right)+D_{t,s}\abs{y}} \]
and
\[ \int_{\reals^{2d}}h\left(x,y\right)f_{t,s}\left(x,y\right)e^{A_{t,s}\abs{x}+B_{t,s}\left(\abs{x}+\abs{y}+1\right)+D_{t,s}\abs{y}}\ud x\ud y<\infty. \]
We therefore have dominated convergence in~\eqref{example_condition_2}, which completes the argument.

\section{Concluding Remarks}  \label{sec:concluding}
 We first provide an intuitive explanation of why the restriction $H\!>\!\nicefrac{1}{2}$ was necessary.
 As $H$ decreases, both the integrand and the integrator in $\int_0^tF(B_\tau)\ud B_\tau$ become rougher,
 whereas the stochastic integral's existence requires a minimal amount of ``cumulative" regularity among
 both terms (this can be best seen in Z\"{a}hle's definition \cite{Zahle} of the integral).
 The threshold turns out to occur at $H=\nicefrac{1}{2}$ (Brownian motion). When $H<\nicefrac{1}{2}$
 the process $B_t$ itself does not necessarily belong to the domain of the divergence operator $\delta$, nor is it even
 an $\mathcal{H}$-valued random variable. In~\cite{DivergenceExtension} Cheridito and Nualart have extended
 the divergence operator to a larger domain, which includes $B_t$ itself and some functions of it.
 However, many of the formulae and theorems don't carry through, and those which do take a different form.
 The case $H<\frac{1}{2}$ is therefore significantly different.
 \medskip

 \noindent
 Secondly, and as $H$ decreases to $\nicefrac{1}{2}$, one expects the covariance of the fractional stochastic integrals
 to converge to that of the Brownian stochastic integrals, and we now proceed to check, skipping some details, that this indeed turns out to be the case, namely, that the right hand side of~\eqref{form} converges to
 $\int_0^{s\wedge t}E_{\nicefrac{1}{2}}F(B_\tau)G(B_\tau)\,\ud\tau$. (Here and henceforth in this section, we add the underlying Hurst parameter as a subscript wherever relevant.) Consider first the case\ $F(x)\!=\!G(x)\!=\!\sgn{x}$\ (for $d\!=\!1$) studied in Section~\ref{Bessel}, for which we now show that the expression in~\eqref{sgncovar} -whose terms are defined in~\eqref{alpha}, \eqref{gamma} and \eqref{eq_sgn.der}-\eqref{corell}\ -\ indeed converges to $s\wedge t$ as $H\!\to\nicefrac{1}{2}^+$.
  \begin{itemize}
  \item[$i$)] $\displaystyle{\lim_{H\to\nicefrac{1}{2}^+}\
                                      {\mathcal P}_H|\tau\!-\!\sigma|^{2H\!-\!2}=\delta(\tau\!-\!\sigma)}$
      in the sense of distributions, and it follows that the first term in the right hand side of~\eqref{sgncovar} converges to\ $s\wedge t$.
  \item[$ii$)] $\displaystyle{\lim_{H\to\nicefrac{1}{2}^+}{\mathcal P}_H\gamma_H(\tau,\sigma)=0}$\ \ a.e.
  \item[$iii$)] $\displaystyle{\lim_{H\to\nicefrac{1}{2}^+}P_H(\tau,\sigma)
                     =\frac{1}{\pi\sqrt{\tau\wedge\sigma(\tau\vee\sigma\!-\!\tau\wedge\sigma)}}}$.
  \end{itemize}
   The last two limits imply that the second integrand in~\eqref{sgncovar} converges to $0$ a.e. It remains to find an integrable upper bound for\ \ $\frac{1}{\tau^H\sigma^H\sqrt{1\!-\!\rho^2_H(\tau,\sigma)}}$,\ uniformly for all $H\!\in\![\nicefrac{1}{2},\nicefrac{3}{4}]$, in order to conclude, by dominated convergence, that the integral of the second term converges to $0$.
  \begin{itemize}
  \item[$iv$)] Note that ${\mathcal P}_H\gamma_H(\tau,\sigma)$ is globally bounded. Passing from $(\tau,\sigma)$ to polar coordinates $(r,\theta)$, and splitting the domain $\tau\!>\!\sigma\!>\!0$ into the two regions $\{\theta\!\in\!(0,\frac{\pi}{6})\}$ and $\{\theta\!\in\![\frac{\pi}{6},\frac{\pi}{4})\}$, one obtains the estimate
 \begin{equation*}
  \frac{1}{\tau^H\sigma^H\sqrt{1-\rho^2_H\left(\tau,\sigma\right)}}\leq
     \ C\left(r^{-1}\!\vee\!r^{-\frac{3}{2}}\right)\ \left\{\begin{array}{ll}
         \theta^{-\frac{3}{8}} & \ \ \ {\rm if}\ \theta\!\in\!(0,\frac{\pi}{6}) \\ \\
         \left(\frac{\pi}{4}-\theta\right)^{-\frac{3}{4}} & \ \ \ {\rm if}\
         \theta\!\in\!(\frac{\pi}{6},\frac{\pi}{4})
                                        \end{array}\ \ .  \right.
 \end{equation*}
  Multiplied by the Jacobian $r$, this bound is in $L^1\!\left((0,R)\mbox{\scriptsize $\times$}(0,\frac{\pi}{4})\right)$ for any $R$.
 \end{itemize}
%
 \bigskip

 \noindent
   The general case, when $F$ and $G$ are arbitrary functions which satisfy the assumptions of Theorem~\ref{generalization}, follows from this particular example after realizing that there exist two positive constants $C_1$ and $C_2$ such that
 \begin{eqnarray*}
  \left|F'\left(\ud x\right)\otimes G'\left(\ud y\right)\left(f_{\tau,\sigma}\right)\right|
        &\leq& C_1\norm{f_{\tau,\sigma}}_{\infty}+C_2\\
        &=&\frac{C_1}{4\pi\tau^H\sigma^H\sqrt{1-\rho^2_H\left(\tau,\sigma\right)}}+C_2\ .
 \end{eqnarray*}
\bibliography{bib_covariance}

\def\cprime{$'$}
\begin{thebibliography}{10}

\bibitem{Adler_monograph}
Robert~J. Adler.
\newblock {\em An introduction to continuity, extrema, and related topics for
  general {G}aussian processes}.
\newblock Institute of Mathematical Statistics Lecture Notes---Monograph
  Series, 12. Institute of Mathematical Statistics, Hayward, CA, 1990.

\bibitem{DivergenceExtension}
Patrick Cheridito and David Nualart.
\newblock Stochastic integral of divergence type with respect to fractional
  {B}rownian motion with {H}urst parameter {$H\in(0,{1\over2})$}.
\newblock {\em Ann. Inst. H. Poincar\'e Probab. Statist.}, 41(6):1049--1081,
  2005.

\bibitem{Folland}
Gerald~B. Folland.
\newblock {\em Real analysis}.
\newblock Pure and Applied Mathematics (New York). John Wiley \& Sons Inc., New
  York, 1984.
\newblock Modern techniques and their applications, A Wiley-Interscience
  Publication.

\bibitem{Guerra&Nualart}
Jo{\~a}o M.~E. Guerra and David Nualart.
\newblock The {$1/H$}-variation of the divergence integral with respect to the
  fractional {B}rownian motion for {$H>1/2$} and fractional {B}essel processes.
\newblock {\em Stochastic Process. Appl.}, 115(1):91--115, 2005.

\bibitem{Hu_Nualart}
Y.~Hu and D.~Nualart.
\newblock Some processes associated with fractional {B}essel processes.
\newblock {\em J. Theoret. Probab.}, 18(2):377--397, 2005.

\bibitem{Hu&Oksendal&Salopek}
Yaozhong Hu, Bernt {\O}ksendal, and Donna~Mary Salopek.
\newblock Weighted local time for fractional {B}rownian motion and applications
  to finance.
\newblock {\em Stoch. Anal. Appl.}, 23(1):15--30, 2005.

\bibitem{Kolmogorov}
A.~N. Kolmogorov.
\newblock {\em Selected works. {V}ol. {II}}, volume~26 of {\em Mathematics and
  its Applications (Soviet Series)}.
\newblock Kluwer Academic Publishers Group, Dordrecht, 1992.
\newblock Probability theory and mathematical statistics, With a preface by P.
  S. Aleksandrov, Translated from the Russian by G. Lindquist, Translation
  edited by A. N. Shiryayev [A. N. Shiryaev].

\bibitem{Mandelbrot}
Benoit~B. Mandelbrot and John~W. Van~Ness.
\newblock Fractional {B}rownian motions, fractional noises and applications.
\newblock {\em SIAM Rev.}, 10:422--437, 1968.

\bibitem{Nualart}
David Nualart.
\newblock {\em The {M}alliavin calculus and related topics}.
\newblock Probability and its Applications (New York). Springer-Verlag, Berlin,
  second edition, 2006.

\bibitem{Pipiras&Taqqu}
Vladas Pipiras and Murad~S. Taqqu.
\newblock Are classes of deterministic integrands for fractional {B}rownian
  motion on an interval complete?
\newblock {\em Bernoulli}, 7(6):873--897, 2001.

\bibitem{Zahle}
M.~Z{\"a}hle.
\newblock Integration with respect to fractal functions and stochastic
  calculus. {I}.
\newblock {\em Probab. Theory Related Fields}, 111(3):333--374, 1998.

\end{thebibliography}
\bibliographystyle{plain}

\end{document}